\documentclass[11pt]{article}
\usepackage{mathrsfs}
\usepackage{amsmath}
\usepackage{amssymb}
\usepackage{graphicx}
\usepackage{amsfonts}
\def \la{\lambda}
\newtheorem{theorem}{\scshape \mdseries  Theorem}[section]
\newtheorem{lemma}[theorem]{\scshape \mdseries  Lemma}
\newtheorem{coro}[theorem]{\scshape \mdseries  Corollary}

\newtheorem{conj}[theorem]{\scshape \mdseries  Conjecture}

\begin{document}
\title{\LARGE \sf The Distance Coloring of Graphs}
\author{Lian-Ying Miao$^{1,}$\thanks{E-mail address: miaolianying@cumt.edu.cn. Supported by NSFC (11271365) and Fundamental Research Funds for the Central Universities (2012LWB46)}, \
Yi-Zheng Fan$^{2,}$\thanks{Corresponding author. E-mail address:
fanyz@ahu.edu.cn. Supported by NSFC(11071002, 11371028), Program for New
Century Excellent Talents in University (NCET-10-0001), Key Project of Chinese
Ministry of Education (210091), Specialized Research Fund for the
Doctoral Program of Higher Education (20103401110002), Scientific Research Fund for Fostering
Distinguished Young Scholars of Anhui University (KJJQ1001).
}\\
{\small  \it $^1$Institute of Mathematics, China University of Mining and Technology, Xuzhou 221116, P. R. China} \\
  {\small \it $^2$School of Mathematical Sciences, Anhui University, Hefei 230601, P. R. China}
 }
\date{}
\maketitle

\noindent {\bf Abstract:}
Let $G$ be a connected graph with maximum degree $\Delta \ge 3$.
We investigate the upper bound for the chromatic number $\chi_\gamma(G)$ of the power graph $G^\gamma$.
It was proved that $\chi_\gamma(G) \le\Delta\frac{(\Delta-1)^{\gamma}-1}{\Delta-2}+1=:M+1$ with equality if and only $G$ is a Moore graph.
If $G$ is not a Moore graph, and $G$ holds one of the following conditions: (1) $G$ is non-regular, (2) the girth $g(G) \le 2\gamma-1$, (3) $g(G) \ge 2\gamma+2$, and the connectivity  $\kappa(G) \ge 3$ if $\gamma \ge 3$, $\kappa(G) \ge 4$ but $g(G) >6$ if $\gamma =2$, (4) $\Delta$ is sufficiently large than a given number only depending on $\gamma$, then $\chi_\gamma(G) \le M-1$.
By means of the spectral radius $\la_1(G)$ of the adjacency matrix of $G$, it was shown that $\chi_2(G) \le \la_1(G)^2+1$, with equality holds if and only if $G$ is a star or a Moore graph with diameter $2$ and girth $5$, and $\chi_\gamma(G) < \la_1(G)^\gamma+1$ if $\gamma \ge 3$.

%
%
%
%



\section{Introduction}
Let $G=(V(G), E(G))$ be a graph.
A {\it vertex $k$-coloring} of $G$ is a mapping from $V(G)$ to the set $\{1, 2, \cdots, k\}$
such that any two adjacent vertices are mapped to different integers.
The smallest integer $k$ for which a $k$-coloring exists is called the {\it chromatic number} of $G$, denoted by $\chi(G)$.
The {\it $\gamma$-th power of the graph} $G$, denoted by $G^\gamma$, is a graph on the same vertex set as $G$ such that two vertices are adjacent in $G^\gamma$
if and only if their distance in $G$ is at most $\gamma$.
The {\it $\gamma$-distance $k$-coloring}, also called {\it distance $(\gamma,k)$-coloring}, is a $k$-coloring of the graph $G^\gamma$ (that is,
any two vertices within distance $\gamma$ in $G$ receive different colors).
The {\it $\gamma$-distance chromatic number} of $G$ is exactly the chromatic number of $G^\gamma$, denoted by $\chi_\gamma(G)$.
Clearly $\chi(G)=\chi_1(G)\leq \chi_\gamma(G)=\chi(G^\gamma)$.

The distance coloring was introduced by Florica Kramer and Horst Kramer \cite{kram1,kram2},
and a recent survey on this topic was also given by them; see \cite{kram3} for more details.
The $\gamma$-distance coloring of graphs has a good application in the frequency assignment problem (or radio channel assignment).
Graph coloring formalizes this problem well when the constraint is that
a pair of transceivers within distance $\gamma$ cannot use the same channel due to interference; see \cite{molloy2}.
The first reference to appear on coloring squares of planar graphs was by Wegner \cite{wegner}.
He posed the following conjecture.
\begin{conj}{\em \cite{wegner}}
 Let $G$ be a planar graph with maximum degree $\Delta$. Then
 \[
\chi_2(G)\leq\left\{
\begin{array}{ll}
7, & \hbox{if~~} \Delta=3,\\
\Delta+5, & \hbox{if~~} 4\leq\Delta\leq 7,\\
\lfloor\frac{3}{2}\Delta\rfloor+1, & \hbox{if~~} \Delta \geq 8.
\end{array}
\right.
\]
\end{conj}

Some work has been done on the case $\Delta = 3$, as listed in \cite[Problem 2.18]{jensen}.
The conjecture is still open, even for the case of $\Delta= 3$.
Many upper bounds on $\chi_2(G)$ for planar graphs $G$ in terms of $\Delta$ have been obtained in the last about two decades.
The best known upper bound so far has been found by Molloy and Salavatipour \cite{molloy}: $\chi_2(G)\leq\lceil\frac{5}{3}\Delta\rceil+78$.
 Havet et al. \cite{havet} proved that $\chi_2(G)\leq\frac{3}{2}\Delta(1+o(1))$ when $\Delta \to \infty$.

We review some results below on the general $\gamma$-distance chromatic number.
It was noted by Skupie\'n \cite{sku} that the
well-known Brooks' theorem can provide the following upper bound:
$$\chi_\gamma(G)\leq 1 +\Delta(G^\gamma) \le 1+\Delta\sum_{k=1}^{\gamma}(\Delta-1)^{k-1}=1+\Delta\frac{(\Delta-1)^{\gamma}-1}{\Delta-2} \;(\Delta \ge 3).\eqno(1.1)$$
If $G$ is planar, Jendrol and Skupie\'n \cite{jend} improved it as
$\chi_\gamma(G)\leq 6+\frac{3\Gamma+3}{\Gamma-2}((\Gamma-1)^{\gamma-1}-1),$
where $\Gamma=\max\{8,\Delta\}$.
Agnarsson and Halld\'orsson \cite{agn} proved that $G^\gamma$ is $O(\Delta^{\lceil \gamma/2 \rceil})$-colorable for any fixed $\gamma$.
Some authors studied $\chi_\gamma(G)$ for special graphs arisen from applications, such as square lattice \cite{fert} and hexagonal lattice \cite{molloy2}.
Other related work could be found in \cite{chen} and \cite{jagger}.
It was proved by Sharp \cite{sharp} that for fixed $\gamma \ge 2$ the distance coloring problem is polynomial time for
$k \le \lceil 3\gamma/2 \rceil$ and NP-hard for $k > \lceil 3\gamma/2 \rceil$.

Though for some classes of graphs $G$, $\chi_\gamma(G)$ is much smaller relative to the order $\Delta^\gamma$, e.g. planar graphs.
There are a lot of graphs whose distance chromatic numbers have the order $\Delta^\gamma$; e.g. Moore graphs (see Theorem \ref{main1} below).
In this paper we discuss the upper bound of $\chi_\gamma(G)$ for a
general $\gamma\ge 2$. These bounds are investigated in two aspects:
one is to use the maximum degree subject to the conditions such as minimum degree,
girth, connectivity, the other is to use the spectral
radius of the adjacency matrix of $G$.

Some notations are introduced as follows.
Let $G$ be a graph. Denote by $G[U]$ the subgraph of $G$ induced by the vertices of $U \subseteq V(G)$.
 For a vertex $v \in V(G)$, denote by $N_G^k(v)$ the set of vertices of $G$ with distance $k$ from $v$.
 Clearly, $N_G^1(v)=:N_G(v)$ is exactly the {\it neighborhood} of $v$ in $G$.
The {\it degree} of a vertex $v \in V(G)$ is denoted by $d_G(v)$.
The {\it maximum} (resp. {\it minimum}) degree of $G$ is denoted by $\Delta(G)$ (resp. $\delta(G)$).
The {\it distance} between two vertices $u$ and $v$ in $G$ is denoted by $dist_G(u, v)$.
 The {\it diameter}, {\it girth}, {\it connectivity} and {\it clique number} of $G$ are denoted by $diam(G),g(G),\kappa(G),\omega(G)$ respectively.
 If $G$ contains no cycles, then we define $g(G)=+\infty$.
The superscript or subscript $G$ may be omitted when it is non-ambiguous.

\section{The upper bound of $\chi_\gamma(G)$ in terms of maximum degree}

When $\Delta=2$, there exist only two connected graphs of order $n$: the path $P_n$ and the cycle $C_n$.
It is easy to get  the  result:

(1) $\chi_\gamma(P_n)=\min\{n, \gamma+1\}$;

(2) $\chi_\gamma(C_n)=\gamma+1$ if $n \!\!\! \mod (\gamma+1)=0$, and $\chi_\gamma(C_n)=\min\{i+1 \ge \gamma+2|n \!\!\! \mod i \le n/i\}$.

%
%

From now on it is assumed that $\gamma \ge 2$ when discussing $\chi_\gamma(G)$, and all graphs are connected with $\Delta \ge 3$.
We denote by $M$ the maximum possible degree of the graph $G^\gamma$ as follows:
$$M=\Delta\frac{(\Delta-1)^{\gamma}-1}{\Delta-2}.$$
First we characterize the graph $G$ for which $\chi_\gamma(G)$ attains the maximum value $M+1$.

\begin{theorem}\label{main1} Let $G$ be a graph.  Then
$$\chi_\gamma(G) \le M+1,$$
 with equality if and only if $G$ is $\Delta$-regular of order $M+1$, $g(G)=2\gamma+1$, $diam(G) = \gamma$,
 i.e., $G$ is a Moore graph with degree $\Delta$ and diameter $\gamma$.
\end{theorem}

 {\bf Proof.} The upper bound follows from (1.1). We now discuss the equality case.
 For the necessity, assume $\chi_\gamma(G)=M+1\ge \Delta(G^\gamma)+1$.
 Then, by Brooks' Theorem, $G^\gamma$ is a complete graph of order $M+1$, and every vertex is at distance at most $\gamma$ of the other $M$ vertices,
 which is maximum obtained when considering that $G$ is a complete $M$-ary tree.
 Thus for every vertex $v$ in $G$, the graph $G$ is a complete $M$-ary tree except for edges between vertices in $N^\gamma(v)$.
 It follows trivially that $G$ is $\Delta$-regular, $g(G)=2\gamma+1$, $\mbox{diam}(G) = \gamma$, i.e. $G$ is a Moore graph.
 The sufficiency is easily obtained from the fact that $G^\gamma$ is a complete graphs of order $M+1$
 as $\hbox{diam}(G) = \gamma$.\hfill $\blacksquare$

\begin{coro} \label{ka3}
 Let $G$ be a graph with $g(G) \ne 2\gamma+1$. Then $\chi_\gamma(G)\leq M$.
\end{coro}

 {\bf Proof.}
 By Theorem \ref{main1}, $\chi_\gamma(G)\leq M+1$, with equality only if $g(G)=2\gamma+1$.
 So, if $g(G) \ne 2\gamma+1$, surely $\chi_\gamma(G)\leq M$.
%
\hfill $\blacksquare$

Next we give some sufficient conditions for a graph $G$ such that $\chi_\gamma(G) \le M-1$.
We will use the idea of {\it saving a color} at a vertex $v$, motivated by Cranston and Kim's work \cite{cra} on list-coloring the square of
subcubic graphs.
A {\it partial (proper) coloring} is the same as a proper coloring except that some vertices may be uncolored.
Given a graph $G$ and partial coloring of $G^\gamma$, we define
$excess(v)$ to be $1+\mbox{(the number of colors available at vertex $v$)}-\mbox{(the number of uncolored neighbors of $v$ in $G^\gamma$)}.$
Note that for any graph $G$ and any such partial $(M-1)$-coloring, every vertex $v$ has $excess(v) \ge 0$.
Similar to the proof of Lemmas 3 and 4 in \cite{cra}, we easily get the following two lemmas.

\begin{lemma}\label{subgraph}
For any edge $uv$ of $G$, $\chi(G^\gamma-\{u,v\}) \le M-1$.
\end{lemma}

\begin{lemma}\label{mainLemma}
Let $G$ be a graph with a partial $(M-1)$-coloring.
Suppose that $u$ and $v$ are uncolored, are adjacent in $G^\gamma$, and that $excess(u) \ge 1$ and $excess(v) \ge 2$.
If we can order the uncolored vertices so that each vertex except $u$ and $v$ is succeeded in the order by at least two adjacent vertices in $G^\gamma$,
then we can finish the partial coloring.
\end{lemma}


%


A simple but useful instance where Lemma \ref{mainLemma} applies is when the uncolored vertices induce a connected subgraph and the vertices $u$ and $v$
are adjacent in the subgraph. In this case, we order the vertices by decreasing distance (within the subgraph) from the edge $uv$.

\begin{coro}\label{mindegree}
 If $G$ is a non-regular graph, then $\chi_\gamma(G)\leq M-1$.
 \end{coro}

{\bf Proof.} Let $v$ be a vertex with minimum degree $\delta \le \Delta-1$, and
let $u$ be a neighbor of $v$.
Note that $d_{G^\gamma}(u) \le M-[1+(\Delta-1)+\cdots+(\Delta-1)^{\gamma-2}] \le M-1$ and
$d_{G^\gamma}(v) \le M-[1+(\Delta-1)+\cdots+(\Delta-1)^{\gamma-1}] \le M-\Delta \le M-3$.
So, if giving a partial $(M-1)$-coloring of $G^\gamma-\{u,v\}$, then $excess(u) \ge 1$ and $excess(v) \ge 2$.
The result follows by Lemma \ref{mainLemma}. \hfill $\blacksquare$

%

\begin{coro} \label{g2d}
If $G$ is a graph with $g(G)\leq 2\gamma-1$, then $\chi_\gamma(G)\leq M-1$.
\end{coro}

{\bf Proof.}
Let $C$ be a cycle of $G$ with length $g(G)$.
For any two adjacent vertices $u$ and $v$ on $C$, $d_{G^\gamma}(u) \le M-2$ and $d_{G^\gamma}(v) \le M-2$.
So, if giving a partial $(M-1)$-coloring of $G^\gamma-\{u,v\}$, then $excess(u) \ge 2$ and $excess(v) \ge 2$. The result follows by Lemma \ref{mainLemma}.
\hfill $\blacksquare$

 \begin{coro} \label{ka5}
 Let $G$ be a graph with $g(G)\geq 2\gamma+2$.
 If $\gamma \ge 3$ and $\kappa(G)\geq 3$, or $\gamma=2$ and $\kappa(G)\geq 4$ but $g(G) > 6$, then $\chi_\gamma(G)\leq M-1$.
\end{coro}

 {\bf Proof.}
 First suppose that $\gamma \ge 3$ and $\kappa(G)\geq 3$.
 Let $C$ be a cycle of $G$ with length  $g(G)$.
 Arbitrarily choose a path $P$ of length $\gamma$ on $C$ connecting two vertices $v_1$ and $x_1$, and choose a vertex $v_2$ lying on $P$ that is adjacent to $v_1$.
 Let $y_1,y_2$ be two vertices of $C$ outside $P$, where $y_1$ is adjacent to $v_1$ and $y_2$ is adjacent to $y_1$.
 As $\kappa(G)\geq 3$ and $g(G) \ge 2\gamma+2$, there exists a path $P'$ of length $\gamma-1$ outside $C$  connecting $v_1$ and $x_2$.
 See the left graph in Fig. 2.1 for these labeled vertices.

We color $x_1, y_1$ with color $1$, and $x_2,y_2$ with color $2$.
As $\kappa(G)\geq 3$, the subgraph $G-\{y_1,y_2\}$ is connected.
Order the vertices of $G-\{y_1,y_2\}$ by decreasing distance (within the subgraph) from the edge $v_1v_2$.
Then $x_1,x_2$ are in the distance class $\gamma-1 \ge 2$.
For any uncolored vertex $w$ other than $v_1$ and $v_2$,
$w$ is always succeeded by at least two adjacent (uncolored) vertices in $G^\gamma$.
Since $v_2$ has two neighbors $x_2,y_2$ in $G^d$ colored by the same color, we have $excess(v_2) \ge 1$.
Since $v_1$ has $4$ neighbors $x_1,y_1,x_2,y_2$ in $G^d$ colored by $2$ colors, we have $excess(v_1) \ge 2$.
The result follows by Lemma \ref{mainLemma}.

Next suppose $\gamma=2$ and $\kappa(G)\geq 4$ but $g(G) > 6$.
By a similar discussion as the above, $G$ contains an induced subgraph $G_2$ as listed in the right side of Fig. 2.1.
We color $x_1,x_2,x_3$ with color $1$.
As $\kappa(G)\geq 4$, the subgraph $G-\{x_1,x_2,x_3\}$ is connected.
Order the vertices of $G-\{x_1,x_2,x_3\}$ by decreasing distance (within the subgraph) from the edge $v_1v_2$.
Then every vertex except $v_1,v_2$ is succeeded by at least two adjacent (uncolored) vertices in $G^\gamma$.
Note $excess(v_2) \ge 1$ and $excess(v_1) \ge 2$.
We can greedily finish the coloring by Lemma \ref{mainLemma}.
%
%
%
%
%
 \hfill $\blacksquare$

\begin{center}
\includegraphics[scale=.6]{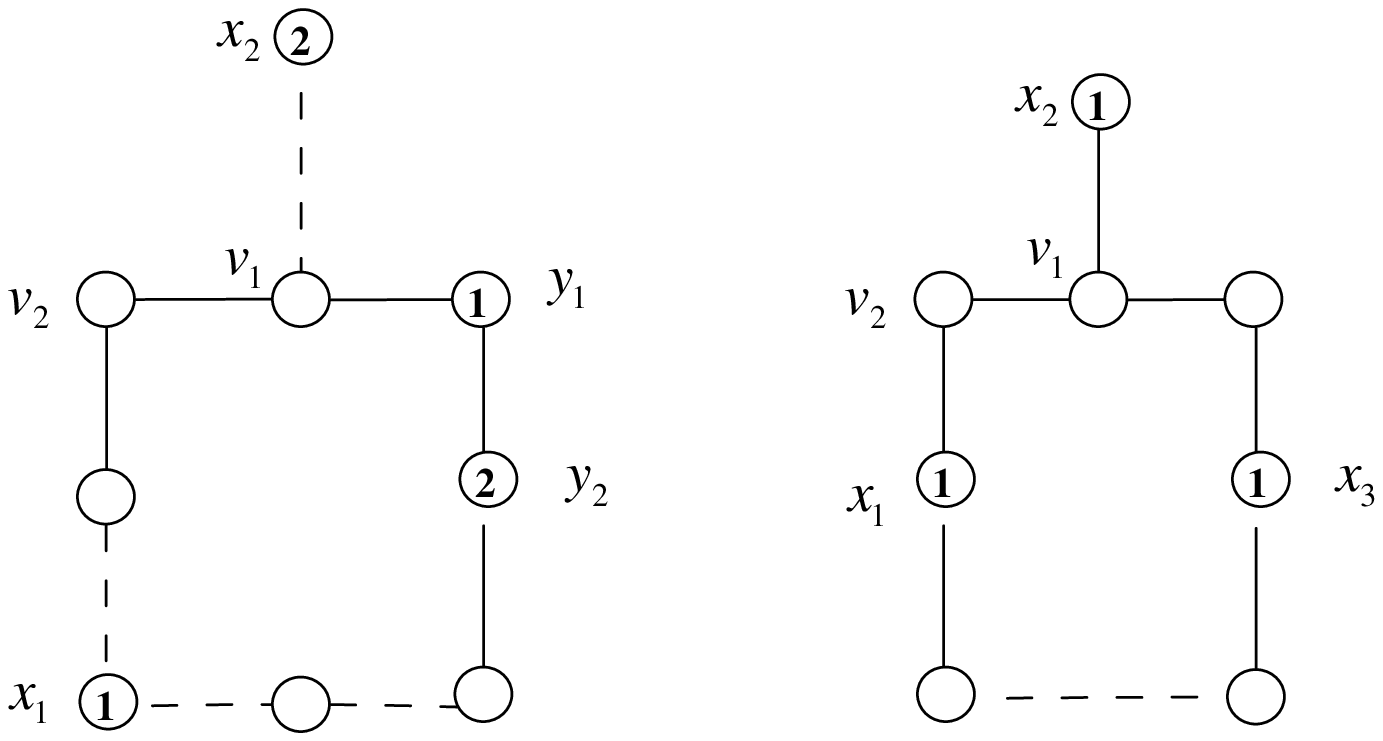}

\vspace{2mm}
{\small Fig. 2.1 The graph $G_1$ (left side) with $dist(v_1,x_1)=\gamma$ and $dist(v_1,x_2)=\gamma-1$, and the graph $G_2$ (right side)}
\end{center}

\begin{lemma} {\em\cite{reed}} \label{reed}
There is a $\Delta_0$ such that if $\Delta(G) \ge \Delta_0$ and $\omega(G) \le \Delta(G)-1$ then
$\chi(G) \le \Delta(G)-1$. In fact, $\Delta_0 = 10^{14}$ will do.
\end{lemma}

\begin{lemma}\label{clique}
If $G$ is a connected graph with maximum degree $\Delta \ge 3$ and $G$ is not a Moore graph, then $G^\gamma$ cannot properly contain a clique of size $M$.
\end{lemma}

 {\bf Proof.}
 If $G$ is a counterexample, that is, $G^\gamma$ properly contains a copy of $K_M$, denoted by $H$.
 Note that $\chi_\gamma(G) \ge M$, so by Corollary \ref{mindegree}, $G$ is a  regular graph.
 Choose adjacent vertices $u$ and $v_1$ such that $v_1 \in V(H)$ and $u \notin V(H)$.
 Note that $|\cup_{j=1}^{\gamma}N^j(v_1)| \le M$, so all vertices in $\cup_{j=1}^{\gamma}N^j(v_1)\setminus \{u\} \subseteq V(H)$.
 Label the neighbors of $u$ as $v_1,v_2,\ldots,v_\Delta$, all belonging to $H$.
 Note that $v_i$ cannot lies on a cycle of length smaller than $2\gamma+1$;
 otherwise, $d_{G^\gamma}(v_i) \le M-1$ and hence $d_{H}(v_i) \le M-2$ as $u \notin V(H)$, which yields a contradiction.

 Now $\cup_{j=1}^{\gamma-1}N^j(v_i)\setminus \{u\} \subseteq V(H)$ for all $v_i$s.
 Since no $v_i$ lies on a cycle of length smaller than $2\gamma+1$, for any $i \ne k$,
 $\cup_{j=1}^{\gamma-1}N^j(v_i)$ and $\cup_{j=1}^{\gamma-1}N^j(v_k)$ share no common vertices except $u$.
 So, the vertices $\left(\bigcup_{i=1}^\Delta \{v_i\} \cup \left(\cup_{j=1}^{\gamma-1}N^j(v_i)\right)\right)\setminus \{u\}$ are exactly $M$ vertices of $H$.
 But $u$ is within distance $\gamma$ of each of these $M$ vertices of $H$.
 Hence, adding $u$ to $H$ yields a clique of size $M+1$, which implies $G$ is a Moore graph.
\hfill $\blacksquare$

\begin{theorem} \label{notexit}
There exists a $\Delta_0$ such that if $\Delta$ is odd and $\Delta \ge \Delta_0$ then
$\chi_\gamma(G) \le M-1$ except $G$ is a Moore graph. In fact, $\Delta_0=(10^{14}+1)^{\frac{1}{\gamma}}+1$ will do.
\end{theorem}

 {\bf Proof.} Assuming to the contrary,
  let $G$ be a non-Moore graph with maximum degree $\Delta \ge \Delta_0=(10^{14}+1)^{\frac{1}{\gamma}}+1$ and $\chi_\gamma(G) \ge M$.
  Note $M$ is an odd number as $\Delta$ is odd, and $M >  (\Delta-1)^{\gamma}-1 \ge 10^{14}$.
  By Corollary \ref{mindegree}, $G$ is $\Delta$-regular.

  By Brooks' theorem, $M \ge \Delta(G^\gamma)\ge \chi(G^\gamma)-1 \ge M-1$.
  So we have three cases: (1) $\Delta(G^\gamma)=M-1$ and  $\chi(G^\gamma)=M$, (2) $\Delta(G^\gamma)=M$ and $\chi(G^\gamma)=M$, and (3) $\Delta(G^\gamma)=M$ and $\chi(G^\gamma)=M+1$.
  For the case (1), also by Brooks' theorem, $G^\gamma$ is a complete graph of order $M$.
  But the sum of the degrees of vertices of $G$ is $M \cdot \Delta$, which is an odd number; a contradiction.
  For the case (2), by Lemma \ref{reed}, $G^\gamma$ properly contains a clique of size $M$, which contradicts to Lemma \ref{clique}.
  The last case cannot occur; otherwise by Theorem \ref{main1} $G$ is a Moore graph.
  \hfill $\blacksquare$

\vspace{3mm}
{\bf Remark 1:}
In this section we give some sufficient conditions for a graph $G$ such that $\chi_\gamma(G) \le M-1$.
We suspect that there exit no graphs $G$ with $\chi_\gamma(G) =M$.

\vspace{3mm}
\noindent
{\scshape \mdseries Conjecture 1}
{\it If $G$ is a connected graph with maximum degree $\Delta \ge 3$ and $G$ is not a Moore graph, then
$\chi_\gamma(G) \le M-1$.}
\vspace{3mm}

Cranston and Kim \cite{cra} conjectured that $\chi_l(G^2) \le \Delta^2-1$, where $\chi_l(G^2)$ is the list-chromatic number of $G^2$.
If their conjure is true, then we will have $\chi_2(G)=\chi(G^2) \le \chi_l(G^2) \le \Delta^2-1$, which implies that Conjecture 1 will hold for $\gamma=2$.

\vspace{3mm}
\noindent
{\scshape \mdseries Conjecture 2}
{\it If $G$ is a connected graph with maximum degree $\Delta \ge 3$ and $g(G)=2\gamma$, then
$\chi_\gamma(G) \le M-1$.}
\vspace{3mm}

By Corollary \ref{mindegree}, it suffices to consider the regular graphs.
Also, if $v$ lies on a $2\gamma$-cycle, then $v$ cannot lie on any other $2\gamma$-cycle because $d_{G^\gamma}(v) \le M-2$ otherwise.
Surely Conjecture 1 implies Conjecture 2 as $g(G)=2\gamma$ implies that $G$ is not a Moore graph.

\vspace{3mm}
\noindent
{\scshape \mdseries Conjecture 3}
{\it If $G$ is a connected graph with maximum degree $\Delta \ge 3$, then $G^\gamma \ne K_M$.}
\vspace{3mm}

Erd\"os et al. \cite{erdos} have proved Conjecture for the case of $\gamma=2$.
If Conjecture 3 is true, then Lemma \ref{clique} would be: $\omega(G^\gamma) \le M-1$, which will generalize the Lemma 22 of \cite{cra};
also we can delete the limitation that $\Delta$ is odd in Theorem \ref{notexit}

If $G$ is a counterexample of Conjecture 3, then it suffices to consider regular graph with girth $g(G) \ge 2\gamma$ by Corollaries \ref{mindegree} and \ref{g2d}.
Also we find $diam(G)=\gamma$.
So $g(G)=2\gamma$, and each vertex lies on exactly one  $2\gamma$-cycle.
However, if we can prove Conjecture 2, then $\chi_\gamma(G) \le M-1$, and hence $G^\gamma \ne K_M$.
So Conjecture 2 implies Conjecture 3.

If $\Delta$ is sufficiently large, then Conjecture 3 also implies Conjecture 1, and hence three conjectures are equivalent.
Let $G$ be a counterexample of Conjecture 1. Then $\chi_\gamma(G) = M$, which implies $G$ is regular and $g(G) \ge 2 \gamma$ by Corollaries \ref{mindegree} and \ref{g2d}.
Hence, $\Delta(G^\gamma)$ equals $M-1$ or $M$.
If  $\Delta(G^\gamma)=M-1$, then by Brooks' theorem $G^\gamma=K_M$; a contradiction to Conjecture 3.
If $\Delta(G^\gamma)=M$ and $\Delta$ is sufficiently large, then by Lemma \ref{reed}, $G^\gamma$ properly contains $K_M$; a contradiction of Lemma \ref{clique}.

\section{Upper bound of $\chi_\gamma(G)$ in terms of spectral radius}
Let $A$ be an entrywise nonnegative matrix.
By Perron-Frobenius theorem, there exists an entrywise nonnegative eigenvector of $A$ corresponding the largest eigenvalue or spectral radius;
this vector is also called the {\it Perron vector} of $A$.
Let $A=[a_{ij}], B=[b_{ij}]$ be two matrices of same size.
Write $A \le B$ if $a_{ij} \le b_{ij}$ for any $i,j$.

Let $G$ be a graph on $n$ vertices.
The {\it adjacency matrix } of $G$, denoted by $A(G)$, is defined as a symmetric $(0,1)$-matrix of order $n$,
where $A(G)_{uv}=1$ if and only if $u$ is adjacent to $v$.
Observe that  $A(G)^k_{uv}$ is the number of walks with length $k$ in $G$ from $u$ to $v$.
So we have
$$ A(G^\gamma) \le A(G)+A(G)^2+\cdots+A(G)^\gamma. \eqno(3.1)$$
Denote by $\la_1(G)$ (resp. $\la_1(A)$) the spectral radius or the largest eigenvalue of $A(G)$ (resp. a square matrix $A$).
The {\it degree matrix} $D(G)$ of $G$ is a diagonal matrix whose diagonal elements are the degrees of the vertices of $G$.
The {\it Laplacian matrix} $L(G)$ of $G$ is defined as $L(G)=D(G)-A(G)$, which is singular and positive semi-definite.
If $G$ is connected, then $0$ is a simple eigenvalue of $L(G)$ with the all-one vector ${\bf 1}$ as the corresponding eigenvector.

\begin{lemma} {\em \cite{wilf}} \label{wilf}
Let $G$ be a graph. Then $\chi(G) \le \la_1(G)+1$, with equality if and only if $G$ is a complete graph or an odd cycle.
\end{lemma}

By Lemma \ref{wilf} and (3.1), we easily get
$$\chi_\gamma(G)=\chi(G^\gamma) \le \la_1(G)+\la_1(G)^2+\cdots+\la_1(G)^\gamma+1=\frac{\la_1(G)^{\gamma+1}-1}{\la_1(G)-1}.\eqno(3.2)$$

However, the upper bound in (3.2) is two large. We improve it in the following, before we prove some basic facts.

\begin{lemma} \label{eigprod}
Let $A,B$ be two nonzero symmetric nonnegative square matrix of same order.
Then
$$ \la_1(AB) \le \la_1(A)\la_1(B),$$
with equality if and only if $A^2,B^2, AB$ share a common Perron vector.
\end{lemma}

{\bf Proof.}
There exists a unit Perron vector $x$ such that $\la_1(AB)=x^TABx$.
So, by Cauchy-Schwarz inequality
$$ \la_1(AB)=x^TABx \le \|Ax\| \cdot \|Bx\| =(x^TA^2x)^{1/2} \cdot (x^TB^2x)^{1/2} \le \la_1(A)\la_1(B),$$
with equality if and only if $A^2x=\la_1(A)^2x, B^2x=\la_1(B)^2x$, i.e. $A^2,B^2,AB$ share a common Perron vector.
\hfill $\blacksquare$

\begin{lemma} \label{matrix} Let $G$ be a connected graph on at least $3$ vertices. Then

\noindent
{\em(1)} $A(G^2) \le A(G)^2-L(G)$ with equality if and only if $g(G) \ge 5$.

\noindent
{\em(2)} If $\gamma \ge 3$, then
\begin{align*}
A(G^\gamma) & \le A(G^{\gamma-1})A(G)-A(G)(D(G)-I)-L(G),\\
A(G^\gamma) & \le A(G)A(G^{\gamma-1})-(D(G)-I)A(G)-L(G),
\end{align*}
both with equalities if and only if $g(G) \ge 2\gamma+1$.

\end{lemma}

{\bf Proof.}
(1) We know that $A(G^2) \le A(G)+A(G)^2$.
But $A(G)^2$ contains nonzero diagonal entries.
In fact, $A(G)^2_{vv}=d(v)$.
So $A(G^2) \le A(G) +A(G)^2-D(G)=A(G)^2-L(G)$. Surely both sides have zero trace.
If the equality holds, then $A(G)_{uv} =1$ implies $A(G)^2_{uv} =0$, that is, $G$ contains no $C_3$.
Furthermore, if $A(G)^2_{uv} \ne 0$, then $A(G)^2_{uv} =1$, i.e. there is exactly one path of length 2 from $u$ to $v$.
Hence $G$ contains no $C_4$.

Conversely, if $G$ contains no $C_3$ or $C_4$, then for each edge $uv$ of $G^2$ (equivalently $A(G^2)_{uv}=1$), $u$ and $v$
are joined by an edge or by exactly one path of length 2 (but cannot happen at the same time).
So $[A(G)^2-L(G)]_{uv}=1$.
If $uv\;(u \ne v)$ is not an edge of $G^2$ (equivalently $A(G^2)_{uv}=0$), then $d(u,v) \ge 3$, and hence $[A(G)^2-L(G)]_{uv}=0$.

(2) Assume that $uv$ is an edge in $G^\gamma$.
If $uv$ is not an edge in $G$, there exists a path $W$ of length $l$ connecting $u$ and $v$, where $2 \le l \le \gamma$.
Writing $W=uw_1\cdots w_{l-1}v$, then $uw_{l-1}$ is an edge of $G^{\gamma-1}$ and $w_{l-1}v$ is an edge of $G$, and hence
$A(G^\gamma)_{uv} \le [A(G^{\gamma-1})\cdot A(G)]_{uv}$.
Surely, $uw_1$ is an edge of $G$ and $w_1v$ is an edge of $G^{\gamma-1}$, and therefore $A(G^\gamma)_{uv} \le [A(G)\cdot A(G^{\gamma-1})]_{uv}$.
Observing that the $uv$-entries of $A(G)(D(G)-I),(D(G)-I)A(G),L(G)$ are all zeros, so two inequalities hold in this case.
If $uv$ is also an edge of $G$, then
\begin{align*}
[A(G^{\gamma-1})A(G)-A(G)(D(G)-I)]_{uv} & \ge \sum_{w \in N(v)\backslash \{u\} } A(G^{\gamma-1})_{uw} A(G)_{wv}-[d(v)-1]=0;\\
[A(G)A(G^{\gamma-1})-(D(G)-I)A(G)]_{uv} & \ge \sum_{w \in N(u)\backslash \{v\} } A(G)_{uw}A(G^{\gamma-1})_{wv}-[d(u)-1])=0.
\end{align*}
So the two inequalities also hold in this case.
The equality cases can be checked directly. \hfill $\blacksquare$

Let $G$ be a graph and let $v \in V(G)$.
The {\it $2$-degree} of $v$ in $G$, denoted by $d^{(2)}(v)$, is defined by $d^{(2)}(v)=\sum_{u \in N(v)}d(u)$.
One can find $d^{(2)}(v)$ is exactly the the sum of the row of $A(G)^2$ corresponding to the vertex $v$.
The graph $G$ is called {\it $2$-degree regular} if all vertices have the same $2$-degrees, or equivalently $A(G)^2$ has the constant row sums.
If the girth $g(G) \ge 5$, then $d^{(2)}(v)=d_{G^2}(v)$.
So, in this case, $G$ is $2$-degree regular if and only if $G^2$ is regular.

\begin{lemma} \label{spectrum}
Let $G$ be a connected graph on at least $3$ vertices. Then

{\em(1)} $\la_1(G^2) \le  \la_1(G)^2 $ with equality if and only if $G$ is $2$-degree regular with $g(G) \ge 5$.
In particular, if, in addition, $G$ is non-bipartite, then the equality holds if and only if $G$ is regular with $g(G) \ge 5$.

{\em(2)} If $\gamma \ge 3$, $\la_1(G^\gamma) < \la_1(G^{\gamma-1})\la_1(G)$ and $\la_1(G^\gamma) <  \la_1(G)^\gamma$.

\end{lemma}

{\bf Proof.}
(1) By Lemma \ref{matrix}(1) and the theory of nonnegative matrices
$$\la_1(A(G^2)) \le \la_1(A(G)^2-L(G)),\eqno(3.3)$$
with equality if and only if $A(G^2) = A(G)^2-L(G)$ as $G^2$ is connected or $A(G^2)$ is irreducible.
Hence, also by Lemma \ref{matrix}(1), the equality in (3.3) holds if and only if $g(G) \ge 5$.

By Wely's inequality,
$$ \la_1(A(G)^2-L(G)) \le \la_1(A(G)^2)+\la_1(-L(G)) = \la_1(G)^2,\eqno(3.4)$$
with equality if and only if $A(G)^2-L(G)$, $A(G)^2$, and $-L(G)$ share a common eigenvector corresponding to their largest eigenvalues.
But $L(G)$ has a simple least eigenvalue $0$ with ${\bf 1}$ as the corresponding eigenvector.
 So the equality in (3.4) holds if and only if $A(G)^2$ has ${\bf 1}$ as an eigenvector, or equivalently $G$ is $2$-degree regular.
The frist part of assertion (1) now follows by combining the above discussion.

If $G$ is non-bipartite, then $A(G)^2$ is irreducible, and has a unique Perron vector (up to multiples), which is necessarily the Perron vector of $A(G)$.
So, in this case the equality  in (3.4) holds if and only if ${\bf 1}$ is a Perron vector of $A(G)$, which implies that $G$ is regular.

(2) By Lemma \ref{matrix}(2),
$$2 A(G^\gamma) \le A(G^{\gamma-1})A(G)+A(G)A(G^{\gamma-1})-A(G)(D(G)-I)-(D(G)-I)A(G)-2L(G).$$
 By Wely's inequality,
\begin{align*}
2\la_1(G^\gamma) & \le \la_1\left(A(G^{\gamma-1})A(G)+A(G)A(G^{\gamma-1})-(A(G)(D(G)-I)+(D(G)-I)A(G))\right)\\
& < \la_1\left(A(G^{\gamma-1})A(G)+A(G)A(G^{\gamma-1})\right)\\
& \le \la_1\left(A(G^{\gamma-1})A(G)\right)+\la_1\left(A(G)A(G^{\gamma-1})\right)\\
& \le 2 \la_1(G^{\gamma-1}) \la_1(G),
\end{align*}
where the last inequality follows from  Lemma \ref{eigprod}.
So, by the first result,
$$\la_1(G^\gamma) < \la_1(G^{\gamma-1})\la_1(G) \le \la_1(G^{\gamma-2})\la_1(G)^2 \le \cdots \le \la_1(G^2)\la_1(G)^{\gamma-2} \le \la_1(G)^\gamma.$$
\hfill $\blacksquare$

\begin{theorem} \label{main2}
Let $G$ be a connected graph on at least $3$ vertices.
Then $\chi_2(G) \le \la_1(G)^2+1,$ with equality holds if and only if $G$ is a star or a Moore graph with diameter $2$ and girth $5$.
If $\gamma \ge 3$, then $\chi_\gamma(G) < \la_1(G)^\gamma+1$.
\end{theorem}

{\bf Proof.}
By Lemmas \ref{wilf} and \ref{spectrum},  $$\chi_d(G)=\chi(G^\gamma) \le \la_1(G^\gamma)+1 \le \la_1(G)^\gamma+1,$$
where the last inequality is strict if $\gamma \ge 3$.
In the case of $\gamma=2$, we have
$$\chi_2(G) \le \la_1(G^2)+1\le \la_1(G)^2+1.\eqno(3.5)$$

If (3.5) holds equalities, then from the first equality, $G^2$ is complete which implies that $diam(G) \le 2$;
and from the second equality,  $G$ is $2$-degree regular and $g(G) \ge 5$.
If $G$ is bipartite, then $G$ must be a tree with diameter at most $2$, which implies $G$ is a star; otherwise, $G$ contains a cycle and $g(G) \ge 6$, which implies $diam(G) \ge 3$, a contradiction.
If $G$ is non-bipartite, then $G$ is regular by Lemma \ref{spectrum}(1).
Hence $G$ is a Moore graph with diameter $2$ and girth $5$.
For the sufficiency, if $G$ is a star, the result holds obviously.
If $G$ is a Moore graph with degree $\Delta$, by Theorem \ref{main1}, $\chi_2(G)=\Delta^2+1$.
Observing $\la_1(G)=\Delta$, we get the equality.
\hfill $\blacksquare$

{\bf Remark 2:}
If $G$ is $\Delta$-regular, then $\chi_2(G) \le \Delta^2+1$ by Theorem \ref{main1}, which is consistent with the bound $\la_1(G)^2+1$ since $\la_1(G)=\Delta$ in this case.
Otherwise, by Corollary \ref{mindegree}, $\chi_2(G) \le \Delta^2-1$.
In this case $\la_1(G) < \Delta$.
A special example is that $G$ is a star on $n$ vertices.
Then Corollary \ref{mindegree} gives $\chi_2(G) \le n^2-2n$, while Theorem \ref{main2} gives $\chi_2(G) \le n$, that latter of which is an equality as the $G^2$ is complete.

\vspace{3mm}
{\bf Acknowledgements:}
We greatly thank an anonymous referee for his invaluable comments, in particular drawing our attention to the reference \cite{cra} from which
we learned a lot of technique to rewrite the proof or strengthen some results in Section 2 of the original version of this paper.

{\small

\end{document}